\documentclass[12pt,reqno]{article}

\usepackage{amsmath}
\usepackage{amsthm}
\usepackage{amscd}
\usepackage{amssymb}
\usepackage{graphicx}
\usepackage{epsf}
\theoremstyle{plain}
\newtheorem{theorem}{Theorem}[section]
\newtheorem{proposition}[theorem]{Proposition}
\newtheorem{lemma}[theorem]{Lemma}
\newtheorem{corollary}[theorem]{Corollary}

\newtheorem{conjecture}[theorem]{Conjecture}

\newtheorem{definition}[theorem]{Definition}
\newtheorem{remark}[theorem]{Remark}

\def\less{<}

\def\int{\text{int}}
\def\invlimit{\smash{\lim\limits_{\raise1pt\hbox{$\longleftarrow$}}}\vphantom{\big(}}
\def\inter{\hskip 1.5pt\raise4pt\hbox{$^\circ$}\kern -1.6ex}
\def\skel(#1,#2){#1^{(#2)}}
\def\hyp {\hbox {\rm {H \kern -2.8ex I}\kern 1.25ex}}
\def\reals {\hbox {\rm {R \kern -2.8ex I}\kern 1.15ex}}
\def\integers {\hbox {\rm { Z \kern -2.8ex Z}\kern 1.15ex}}
\def\naturals {\hbox {\rm {N \kern -2.8ex I}\kern 1.20ex}}
\def\rationals {\hbox {\rm { Q \kern -2.2ex l}\kern 1.15ex}}
\def\hyp {\hbox {\rm {H \kern -2.7ex I}\kern 1.25ex}}

\begin{document}

\title{On the intersection of unknotting tunnels and the decomposing annulus in connected sums}

\author{Yoav Moriah\thanks{Supported by The Fund for Promoting Research at the
Technion,  grant 100-127 and the Technion VRP fund, grant 100-127.}}

\date{}

\maketitle

\begin{abstract}
Given $(V_1,V_2)$  a  Heegaard splitting of  the complement of a composite knot $ K = K_1 \# K_2$  in 
$S^3$, where $K_i, i = 1,2$ are prime  knots, we have  a unique, up to isotopy,  decomposing  annulus $A$. 
When the intersection of  $A$ and $V_1$ is a  minimal collection of disks we study the components of 
$V_1 - N(A)$ and show  that at most one component is a  $3$-ball meeting $A$ in two disks. This is a crucial
step in proving the conjecture that a necessary and sufficient condition for the tunnel number of a connected sum
to be less than or equal to the sum of the tunnel numbers is that one of the knots has a Heegaard splitting in which
a merdian curve is primitive.
\end{abstract}

\section{Introduction}
\label{intro}

\vskip15pt

The way in which the tunnel number $t(K)$ of a knot $K = K_1 \# K_2$   relates to $t(K_1)$ and 
$t(K_2)$  is a long standing question. It was long known that $t(K_1 \# K_2) \leq t(K_1) + t(K_2) +1$. 
That this inequality is best possible was proved by the second author and Rubinstein in  [MR] and 
by  Morimoto, Sakuma and Yokota in [MSY]. However it is not yet undersrtood when this phenomenon 
occurs.  We state the following conjecture and note that it was proved to be true by Morimoto, for the 
special class of knots which do not contain essential meridional surfaces, also known as {\it smallish knots}
(see  [Mo3]): 

\begin{conjecture}\rm  
The knots $K_1 \subset S^3,  K_2  \subset S^3$ and $K_1 \# K_2 $ satisfy the inequality  $t(K_1 \# K_2) \leq
t(K_1) + t(K_2)$ if and only if either $E(K_1)$ or $E(K_2)$ , say  $E(K_1)$, has a minimal genus Heegaard
splitting $(V_1,V_2)$  in which  $\partial E(K_1) \subset V_1$ and a  meridian curve $\mu \subset \partial
E(K_1)$ is isotopic to a curve $\mu^*$ on the Heegaard surface $\partial V_2$ so that $\mu^*$ intersects an
essential disk $D \subset V_2$ in a single point.

\end{conjecture}

A meridian which has the property above will be called {\it primitive}.

\vskip10pt

The \lq\lq  if" part of the statement is a well known result proved in Lemma \ref{onlyif}, however the
\lq\lq only if " part turned out to be much more difficult.

\vskip10pt

When attempting to prove the conjecture one needs to deal with the following stituation :

We are given a minimal genus Heegaard  splitting $(V_1,V_2)$ for $E(K_1 \# K_2)$ and we need to find
\underbar{minimal} genus Heegaard splittings for $E(K_1)$ or $E(K_2)$ in which a meridian is primitive. We
may  assume that $V_1$ is a small neighborhood of the spine i.e., $V_1$ is a regular neighborhood of the
boundary torus union a regular neighborhood of the tunnels, and we can choose a decomposing annulus $A$
which intersects $V_1$ in a  minimal number of disks. However it is exactly here that a major problem arises as it
is concievable, a priori, that the interior of a single tunnel will intersect $A$  in a large number of points. If this
occurs, the decomposing annulus $A$ will cut $\partial V_1$ into pieces, a lot of which will be annuli. Since the
Euler charecteristic of an annulus is zero we can have arbitrarily many of them. It  follows that the induced (as in
Section
\ref{prelim}) Heegaard splittings on $E(K_1)$ and $E(K_2)$ will be of very high genus and it is easy to generate
such Heegaard splittings where the meridian curve $\mu$ is primitive. However these high genus Heegaard
splittings will give us no information about a minimal genus one.

In this paper we show that this situation  cannot happen.
Let us call a closure of a component of $V_1 - A$ which intersects $A$ in $n$ disks an {\it n-float}. With
this terminology we show:

\vskip10pt
\noindent {\bf Theorem \ref{UturnProp}.}  {\it Let $K_i \subset S^3$ be prime knots and $K = K_1 \# K_2$ be 
their connected sum. Let $(V_1, V_2)$  be a Heegaard splitting of  $E(K)$ and $A$ a decomposing
annulus which intersects $V_1$ in a minimal number of disks. Then: $V_1 - N(A)$  has at most one 
component which is a $2$-float of genus zero. }
\vskip 25pt

\noindent For the proof we need the following result which is of independent interest:
\vskip 10pt

\noindent {\bf Theorem \ref{innertunnelthm}.}  {\it Let $K_1, K_2 \subset S^3 $ be prime knots. Then evey Heegaard splitting $(V_1, V_2)$ for 
$E(K) = E(K_1 \# K_2)$ has a spine which contains at least one cycle disjoint from a decomposing
annulus for $K$ which minimizes the number of  intersections with  $V_1$.  }

\vskip 10pt

Note that if the meridian $\mu$ of $K_1$ or $K_2$ is primitive it does not follow that the Heegaard splitting
of $E(K_1 \# K_2)$  induced, as in Lemma \ref{onlyif}, by the Heegaard splittings of $E(K_1)$ and  $E(K_2)$ is
of minimal genus. It was  shown by Morimoto in  [Mo2] that there are prime knots for which 
$t(K_1 \# K_2) \less t(K_1) + t(K_2)$. In further work Kobayashi showed that if one allows the knots to be
non-prime then the degeneration of tunnel number can be arbitrary i.e, for each $n$ there are knots $K_1^n$ and
$K_2^n$ so that $t(K_1^n \#K_2^n) \less t(K_1^n) + t(K_2^n) - n$. These results are characterized by an attempt
to get upper bounds on $t(K)$ in terms  of $t(K_1) + t(K_2)$.

Results in the opposite  direction were obtained by  M. Scharlemann and J. Schultens. They  analyze 
essential annuli  coming from the intersection of a decomposing annulus with a strongly irreducible 
Heegaard splitting to get results generalizing those of Kwong (see [SS] and [Kw]). In particular they 
show that  $t(K) \geq {2\over5}( t(K_1) + t(K_2))$.

\vskip10pt

\noindent  For definitions of the above terminology see Section \ref{prelim}.

\bigskip

\noindent {\bf Acknowledgments:} We would like to thank  Bronek Wajnryb for many conversations and
Ying-Qing Wu for suggesting a somewhat shorter version of the proof of Theorem 4.1. Also thanks to Anna
Klebanov for an argument for Case 1 in the proof of Theorem 4.1.
\bigskip
\eject

\section{Preliminaries}
\label{prelim}

\noindent In this section we define some of the notions and state some of the results needed
for the proof of the main theorem.

Throughout the paper  $K_1$ and $K_2$ will be knots in $S^3$ and  $K = K_1 \# K_2$
will denote  the connected sum of $K_1$ and $K_2$.  The knots $K_i$ will be called 
the {\it summands} of the {\it composite knot} $K$. Let $N()$ denote an open 
regular neighborhood in $S^3$.

Recall that $(S^3,K)$ is obtained by removing  from each space $(S^3,K_i),  i = {1,2},$ a 
small $3$-ball intersecting $K_i$ in a short unknotted arc and gluing the two remaining $3$-balls along 
the $2$-sphere boundary so that the pair of points of $K_1$ on the $2$-sphere are identified with  the 
pair of points of $K_2$. If we denote $S^3 - N(K)$ by $E(K)$ then $E(K)$ is obtained from 
$E(K_i), i = {1,2},$ by identifying a meridional annulus $A_1$ on $\partial E(K_1)$ with a 
meridional annulus $A_2$ on $\partial E(K_2)$. A knot $K \subset S^3$ is {\it prime} if it is not a 
connected sum of two non-trivial knots. The annulus $A_1 = A_2$ will be denoted by $A$ 
and called the{ \it decomposing annulus}. If both knots $K_1, K_2$ are prime then the decomposing 
annulus is unique up to isotopy  

 A {\it tunnel system} for an arbitrary knot $K \subset S^3$ is a collection of 
properly  embedded arcs  $\{t_1, \dots, t_n\}$  in  $S^3 - N(K)$  so that 
$S^3 - N(K \cup t_1 \cup \dots \cup t_n)$ is a handlebody.

Given a tunnel system for a knot $K \subset S^3$ note that the closure of 
$N(K \cup t_1 \cup \dots \cup t_n)$ is always a handlebody denoted by $V_1$ 
and the handlebody $S^3 - N(K \cup t_1 \cup  \dots \cup t_n)$ will be denoted by $V_2$. 
For a given knot $K \subset S^3$ the smallest cardinallity of any tunnel system is called the 
{\it tunnel number} of $K$ and
is denoted by $t(K)$.

A compression body  $V$  is a  compact orientable and connected 3-manifold with a preferred 
boundary component $\partial_+V$  that is obtained from a collar of $\partial_+ V$ by attaching 
2-handles and 3-handles, so that the connected components of  $\partial_- V$ = $\partial V - \partial_+ V$ 
are all distinct from  $S^2$.  The extreme cases, where  $V$  is a handlebody i.e., $\partial_- V = \emptyset$,
or where $V = \partial_+V \times I$, are allowed.  Alternatively we can think of $V$ as obtained from
$(\partial_-V) \times I$ and $3$-balls by attaching $1$-handles to $(\partial_-V) \times \{1\}$ and the $3$-balls. 
An essential annulus
in a  compression body will be called a {\it vertical (or a spanning) annulus } if it has its boundary components
on both of $(\partial_-V)$ and $(\partial_+V)$.

 A Heegaard splitting for a manifold $M$ is a decomposition $M = V_1 \cup V_2$, so that   
and $ V_1 \cap V_2 =  \Sigma$, where $V_i$, $ i = 1,2$ are  compression  bodies and  
$\Sigma = \partial_{+} V_1 =  \partial_{+} V_2$ is the Heegaard  splitting surface.

Given a knot  $K \subset S^3$ a {\it Heegaard splitting } for $E(K)$ is a decomposition of $E(K)$ into
a compression body $V_1$ containing $\partial E(K)$ and a handlebody $S^3 - int(V_1)$. Hence, a  
tunnel system $\{t_1, \dots, t_n\}$  in  $S^3 - N(K)$  for $K$ determines a Heegaard splitting of genus
$n +1$ for $E(K)$. Conversely given a Heegaard splitting $(V_1, V_2)$ for $E(K)$ any minimal complete
disk system ${\cal D} = \{D_1, \dots, D_n\}$ so that $V_1 - N({\cal D}) $ = $T^2 \times I$ determines a 
tunnel system by taking $t_i$ to be the cocore arc $\{0\} \times I$ of $cl(N({\cal D}_i)) = {\cal D}_i \times
I$

Given a Heegaard splitting $(V_1, V_2)$  for $ S^3 - N(K_1 \# K_2)$ we will choose  a 
decomposing annulus $A$  which intersects the compression body $V_1$ in two {\it vertical annuli}  
$A^*_1, A^*_2$ i.e., annuli with one boundary component on $\partial_+$ and one on $\partial_-$, and 
in a \underbar{minimal} collection of disks ${\cal D} = \{D_1, \dots, D_d\}$ . Note also that 
$A$ intersects $V_2$ in a connected incompressible planar surface.

Let $ {\cal E} = \{E_1 , \dots, E_{t(K) + 1}\}$ be a complete meridian disk 
system for $V_2$, chosen to minimize the intersection with $A$. Since $V_2$ is a handlebody 
it is irreducible  and we can assume that no component of  ${\cal E} \cap A$ is a simple closed curve.

When we cut $E(K)$ along a decomposing annulus $A$ any Heegaard splitting $(V_1, V_2)$ of $E(K)$
induces  Heegaard splittings on both of $E(K_1)$ and $E(K_2)$, as follows: Set $V_1^i = (V_1 \cap E(K_i))
\cup_{{\cal D} \cup A_1^* \cup A_2^*} N(A)$; it is a compression body as it is a union of an 
$annulus \times I$ and some 1-handles, glued along the two vertical annuli and a collection of disks. Now set 
$V_2^i = V_2 - N(A)$;  it is a handlebody since the annulus $A$  meets $V_2$ in an incompressible
connected planar surface $P$ which separates $V_2$ into two components each of which is a handlebody. 
Hence the pair $(V_1^i,V_2^i)$ is a Heegaard splitting for $E(K_i)$ and will be referred to as the 
{\it induced Heegaard splitting} of $E(K_i)$.

In the other direction we have:
\vskip10pt
\begin{lemma}
\label{onlyif}

Assume that either $E(K_1)$ or $E(K_2)$ , say  $E(K_1)$, has a minimal genus Heegaard
splitting $(V^1_1,V^1_2)$  in which  $\partial E(K_1) \subset V^1_1$ and a  meridian curve $\mu \subset \partial
E(K_1)$ is isotopic to a curve $\mu^*$ on the Heegaard surface $\partial V^1_2$ so that $\mu^*$ intersects an
essential disk $D \subset V^1_2$ in a single point. Then the knots $K_1 \subset S^3,  K_2  \subset S^3$ and 
$K_1 \# K_2 $ satisfy  $t(K_1 \# K_2) \leq t(K_1) + t(K_2)$.

\end{lemma}
\vskip10pt
\begin{proof} Choose a meridional annulus $A_1$ in $E(K_1)$ so that $A_1 \cap V^1_1$ are two vertical annuli
$A_1^*, A_2^*$ and  $A_1 \cap V^1_2$ is a regular neighborhood of $\mu^*$.  Let $(V_1^2, V_2^2)$ be
any Heegaard splitting of $E(K_2)$. When we glue $E(K_1)$ to $E(K_2)$ by identifying the meridional annulus
$A_1$ with a meridonal annulus $A_2 \subset \partial E(K_2)$ we glue $V_1^1$ to $V_1^2$ along two vertical
annuli, so that the result is clearly a compression body $V_1$ with $\partial_ - V_1 = T^2$. We also need to
glue $V_1^2$  to $V_2^2$ by identifying $A_1$ with an annulus on $V_2^2$. Since $A_1$ intersects an
essential disk $D$ in $V_2^1$ in a single arc the  resulting manifold will be a handlebody $V_2$ of genus
$g(V_2^1) + g(V_2^2) - 1 = t(K_1) + t(K_2)$. So we have obtained a Heegaard splitting $(V_1,V_2)$
of $E(K_1 \# K_2)$ of genus $t(K_1) + t(K_2) +1 $ and hence $t(K_1 \# K_2) \leq t(K_1) + t(K_2) .$

\end{proof}
\vskip10pt

Following the notation of Morimoto (see [Mo1]) we consider now the planar surface $ P = A \cap V_2$.
It has two distinguished boundary components coming from the vertical annuli  $A^*_1, A^*_2$ and 
denoted by $C^*_1, C^*_2$ respectively. There are exactly $d$ other boundary components of $P$ which we denote 
by $C_1, \dots, C_d$. With this notation we have $\partial D_i = C_i$.  The arcs of ${\cal E} \cap A$ are 
contained in $P$ and can be divided into three types: 

\begin{itemize}
\item[1.] An arc $\alpha$ of Type I is an arc connecting two different boundary components of $P$.
\item[2.]  An arc $\alpha$ of Type II is an arc connecting a single boundary component of $P$ 
to itself so that the arc does not separate the boundary components $C^*_1, C^*_2$.
\item[3.] An arc $\alpha$ of Type III is an arc connecting a single boundary component of $P$ 
to itself with the additional property that the arc does separate the boundary components     
$C^*_1,   C^*_2$.
\end{itemize}

Since the annulus $A$ was chosen to minimize the number of disks in $A \cap V_1$ the planar surface $P$ is
incompressible in the handlebody $V_2$. Hence there is a sequence of boundary compressions of $P$ 
along disjoint arcs $\alpha_i$ using sub-disks of $\cal E$ so that the end result is a collection of disks. 
Any such sequence defines an order on the  arcs $\alpha_i$ .

\begin{definition}
\label{darcDef}
\rm Let  $\alpha_i$ be an arc of intersection  of $P \cap {\cal E}$. We call $\alpha$  a  $ d-arc$ if there is some
compression order so that $\alpha_i$ is of type I and there is some component $C$ of $\partial P -
(C^*_1\cup C^*_2)$  which  meets $\alpha_i$ and does not meet any $\alpha_j$ for any $j \less i$. If 
$\alpha_i$ is of type I and connects $C^*_1$ to
$C^*_2$ it is called an $e-arc$.

\end{definition}

Any outermost arc $\alpha_i$ determines a sub-disk $\Delta$ on some $E_i$ where $\partial \Delta =
\alpha_i \cup \beta$ and $\beta$ is an arc on $\partial V_1 = \partial V_2$. When we perform an 
isotopy of type A  i.e., pushing $P$ through $\Delta$ as in [Ja], we produce a band $b$ with core
$\beta$ on $\partial V_1 = \partial V_2$. The following result is proved in  [Mo1] pp. 41- 42, 
and [Oc]: 

\begin{theorem}[ Morimoto]
\label{darcThm}
If the decomposing annulus is chosen to minimize the number of
components of $V_1 \cap A$ and $V_1 \cap A \not= A^*_1 \cup A^*_2$  then in $V_2 \cap E$ = $P \cap E$:
\begin{itemize}
\item[\rm (a)] there are no d-arcs. 
\item[\rm (b)] there are no e-arcs.
\item[\rm (c)] there are no arcs of type II.
\item[\rm (d)] each component $C \subset \partial P$ has an arc $\alpha$ of type III with end points
on $C$. 

\end{itemize}

\end{theorem}

\noindent {\bf Remark:} It follows from the above theorem that every component $C$ of 
$\partial P$ has arcs of type I and each component $C$ of $\partial P - (C_1^* \cup C_2^*)$ has  arcs of type III.
Thus there is an order on the disk components of $A \cap V_1$ starting from $A^*_1$ to $A^*_2$ or vice versa.

\bigskip

\section{Interior tunnels}

Consider now a Heegaard splitting $(V_1 , V_2)$  for $E(K)$ the exterior of $K = K_1 \# K_2$, 
where  $\partial E(K) \subset V_1$ and  the decomposing annulus $A$ meets $V_1$ in disks 
and two vertical annuli.  Since the annulus $A$  meets $V_2$ in a connected planar surface $P$
it separates $V_2$ into two components each of which is a handlebody. We will denote the 
handlebodies $cl(V_2 - A) \cap E(K_i)$ by $V_2^i$ respectively. 
However $V_1 - A $ might have many components.

\begin{definition}
\label{floatdef}
\rm A component of $cl(V_1 - A) $ which is disjoint from $\partial E(K_i)$ and 
intersects $A$ in $n$ disks will be called an  {\it n-float}. An $n$-float is either a 3-ball or a handlebody of
some genus $g$ if its spine is not a tree. In this case we say that the $n$-float is of genus $g$.
\end{definition}

\noindent {\bf Remark:}  Note that  there are always 
exactly two components of $cl(V_1 - A) $ not disjoint  from $\partial E(K_i)$ (one in each of $E(K_1)$ and
$E(K_2)$) and each one is a handlebody  of genus at least one as $V_1$ is a compression body with a $T^2$
boundary. We denote these  special  components by $N_1$ and $N_2$ depending on whether they are
contained in  $E(K_1)$ or $E(K_2)$ respectively.

\bigskip
Consider now $E_i \subset \cal E$ any one of the meridian disks of $V_2$. On $E_i$ we have a 
collection of arcs corresponding to the intersection with the decomposing annulus. These arcs, 
as indicated in Fig. 1, separate $E_i$ into sub-disks where disks on opposite sides of arcs  
are contained in opposite sides of $A$ i.e., in $E(K_1)$ or  $E(K_2)$ respectively. So each sub-disk 
is contained in either $E(K_1)$ or  $E(K_2)$. The boundary of these sub-disks is a collection 
of alternating arcs $\cup(\alpha_i  \cup \beta_i)$ where $\alpha_i$ are arcs on $A$ and $\beta_i$
are arcs on some component of  $cl(V_1 - A)$.
\vskip10pt

\begin{figure}[htm]
\centerline{\epsfysize=2.2truein\epsfbox{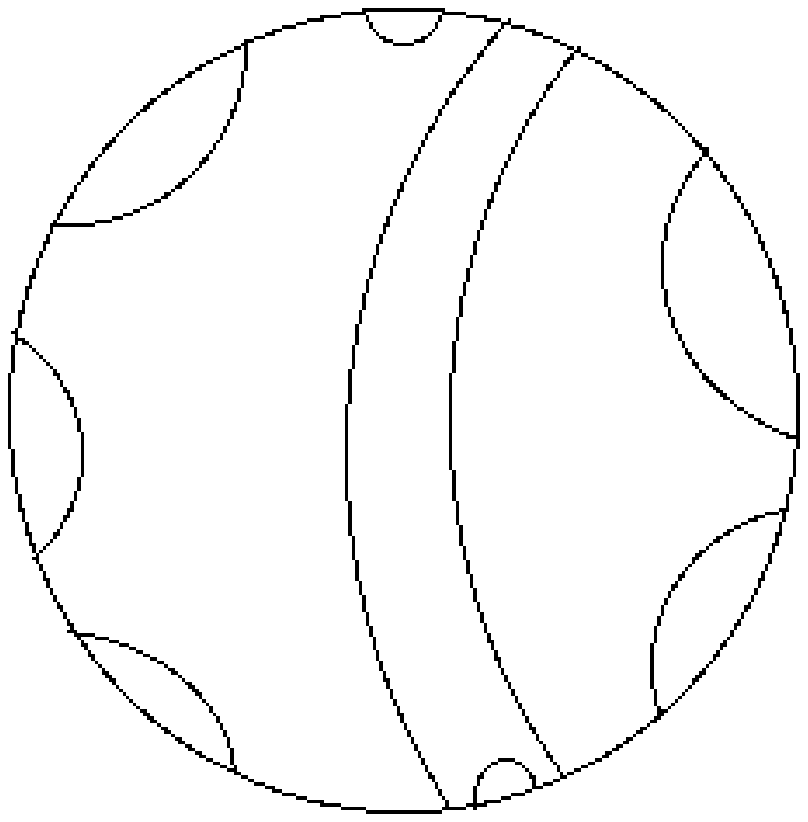}}
\vskip10pt

\centerline{Fig. 1}
\end{figure}
\vskip15pt

\begin{proposition}
\label{innertunnelProp}
Let $K_1$ and $K_2$ be knots in $S^3$ and let $K, A , \cal E$ be the  connected
sum, a minimal intersection decomposing annulus and a meridional system for some Heegaard splitting
of $E(K)$ as above. Then
 
\begin{itemize} 
\item[\rm (a)] the $\beta$ arc part of the boundary of an outermost sub-disk in $E$ cannot
 be contained in a n-float of genus $0$.
\item [\rm (b)]  if  the $\beta$ arc part of the boundary of an outermost 
sub-disk in $E$ is contained in an $N_i$ component $i = 1$ or $2$  and  $K_i$ is prime the
genus of $N_i$ is greater than  one.
\end{itemize}

\end{proposition}

\begin{proof}
Denote an outermost sub-disk of some $E_j$ by $ \Delta$  and suppose it is cut off by an arc $\alpha$ 
on $A$ with end points on a disk $D_i$  which belongs to some n-float of genus $0$. Since $\Delta$ is 
an outermost disk, by Theorem \ref{darcThm}, the arc $\alpha$ must be of type III. Further assume 
$\partial \Delta = \alpha \cup \beta$ where $\beta$ is an arc on the  n-float meeting $D_i$ in exactly two 
points $\partial \beta = \partial \alpha$. On $\partial D_i$  there is a small arc $\gamma$ so that $\gamma 
\cup \beta$ is a simple closed curve  on the n-float bounding a disk $D$ there, since the $n$-float has no
genus (see Fig. 2 below). Furthermore $\gamma \cup \alpha$ is a simple closed loop on $A$ which together 
with a boundary component of $A$ bounds a sub-annulus  of $A$. Hence $\gamma \cup \alpha$  bounds a 
disk $D'$ on the decomposing $2$-sphere of $K$ intersecting $K$ in a single point. Thus we obtain a 
$2$-sphere $D \cup \Delta \cup D'$ which intersects the knot $K$ in a single point. This is a contradiction 
which finishes case (a). 

For case (b), assume that  the outermost disk $\Delta$ is contained in $N_1$, say, and that genus 
$N_1$ is one. As before we have $\partial \Delta = \alpha \cup \beta$ where $\beta$ is an arc on 
$N_1$ and a small arc $\gamma$ so that $ \gamma \cup \beta$ is a simple closed curve on $N_1$. 
If $ \gamma \cup \beta$ bounds a disk in $N_1$ we have the same proof as in case (a). If
 $ \gamma \cup \beta$ does not bound a disk on $N_1$ we consider small sub-arcs $\beta_1$
and $\beta_2$ of $\beta$ which are respective closed neighborhoods of $\partial \beta$. These
arcs together with a small arc $\delta$ on $\partial N_1  - \partial E(K_1)$  and $\gamma$ bound 
a small band $b$ on $\partial N_1$. Notice that $b  \cup_ {\beta_1,\beta_2}\Delta$ is an annulus
$A'$. The annulus $A'$ together with the sub-annulus $A''$ of $A$ cut off by $\alpha \cup \gamma$ 
defines an annulus $A' \cup_{\alpha \cup \gamma} A''$ which determines an isotopy of a meridian 
curve $C_1$ to a simple closed curve $\lambda$ on $\partial N_1$. Note that $N_1$ is a solid torus 
and $\pi_1(N_1) =  \integers$ which is generated by a meridian $\mu$ of $E(K_1)$. Hence 
$[\lambda] =   [C_1] = \mu \in \pi_1(N_1)$ (see Fig. 3).

\begin{figure}[htm]
\hskip60pt{\epsfysize=3.2truein\epsfbox{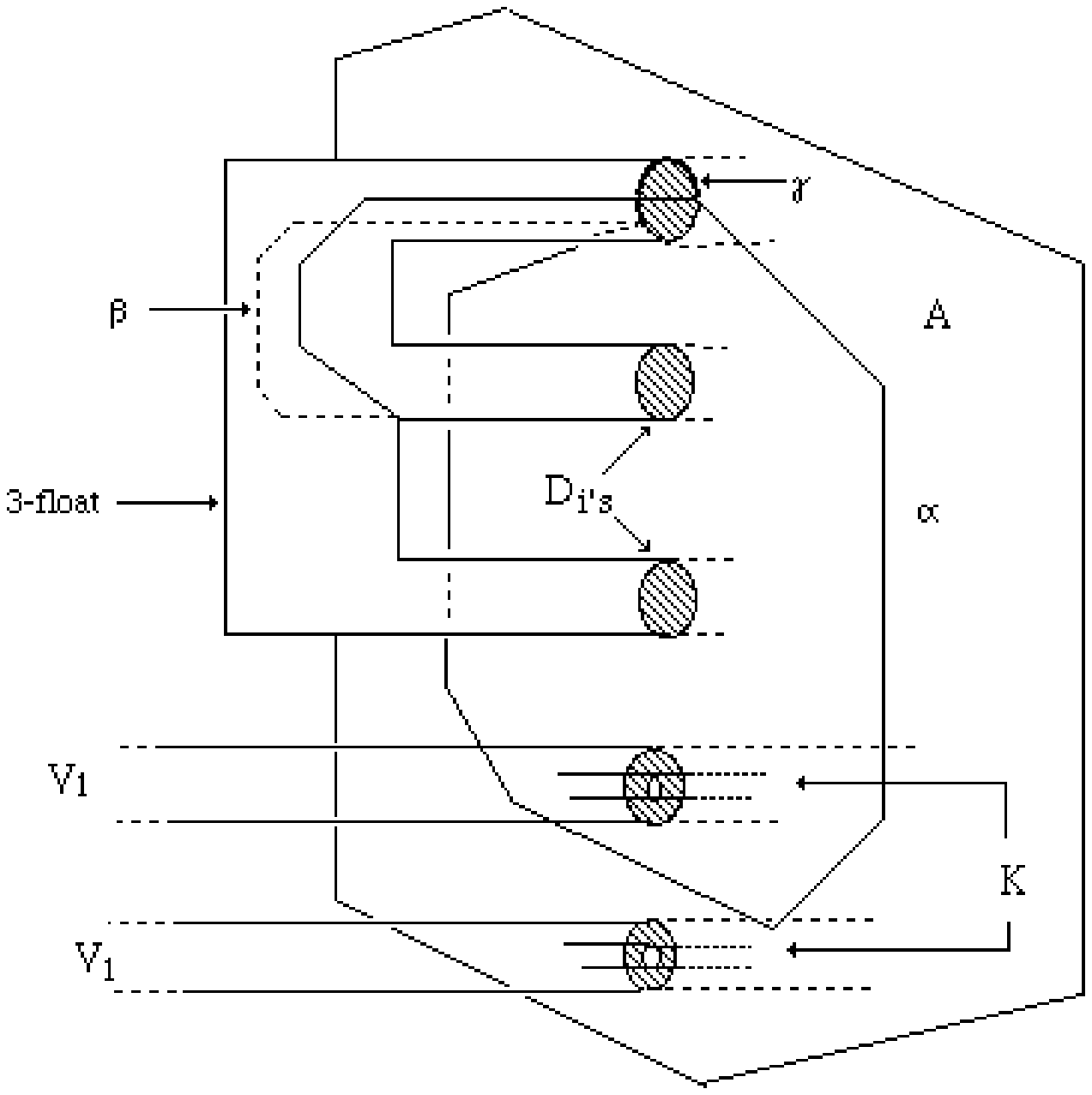}}
\vskip10pt

\centerline {Fig. 2}
\end{figure}

\begin{figure}[htm]

\hskip60pt{\epsfysize=3.2truein\epsfbox{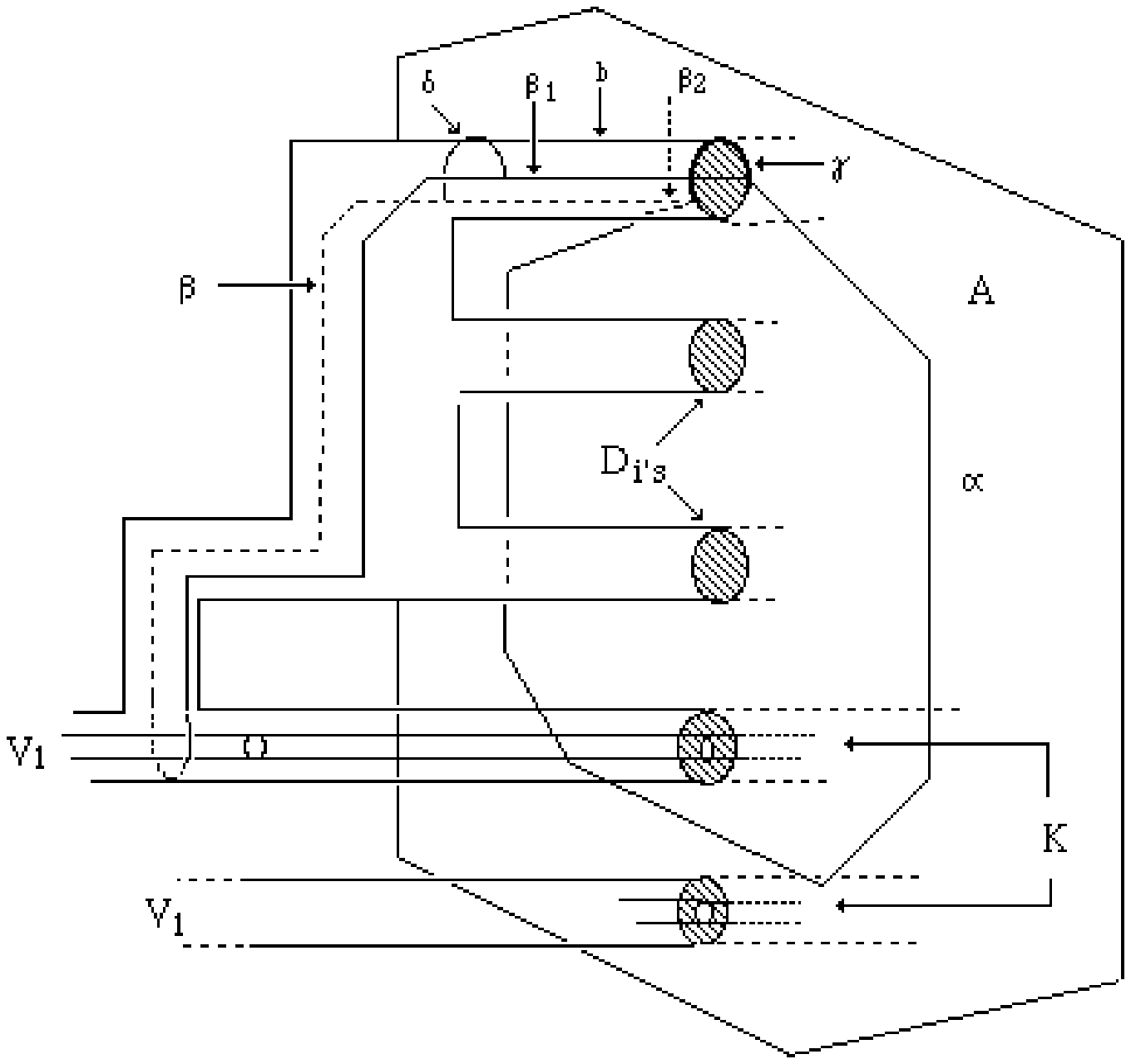}}
\vskip10pt

\centerline {Fig. 3}
\end{figure}

\vskip5pt
Now we can consider the  
annulus $(A - A'') \cup A'$. If it is non-boundary parallel then since both knots $K_1, K_2$ are prime it must be
a  decomposing annulus which has at least one less disk component intersection than $A$ in contradiction to the
choice of $A$.  If it is boundary parallel, then as above, we have $A'' \cup A'$ as a decomposing annulus with a
smaller  number of disks. Again in contradiction  to the choice of $A$. So genus $N_1$ cannot be one and this
finishes case (b).

 \end{proof}
As a corollary we obtain:
\vskip 8pt 
\begin {theorem}

\label{innertunnelthm}

Let $K_1, K_2 \subset S^3 $ be prime knots. Then evey Heegaard splitting $(V_1, V_2)$ for 
$E(K) = E(K_1 \# K_2)$ has a spine which contains at least one cycle disjoint from a decomposing
annulus for $K$ which minimizes the number of  intersections with  $V_1$. 
\end{theorem}

\begin{proof}
Since the $\beta$ part of an outer-most disk must be contained in a float of genus greater than 
one we must have a $1-handle$ on the float to create the genus. The core arc of the $1$-handle 
which is disjoint from the decomposing annulus $A$ determines the cycle.

\vskip10pt 
\end{proof}

\vskip 10pt

\noindent As a side benefit of these considerstions we have the following:
\vskip 10pt

\begin{remark}
\label{benefit} \rm
 
Given a  minimal genus Heegaard splitting $(V_1,V_2)$ for $E(K) = E(K_1 \# K_2)$, where $K_i , i = 1,2$  
are prime knots in $S^3$, assume that the decomposing annulus with minimal intersection with 
the compression body $V_1$  has $d$ disk components and two vertical annuli. Recall that
if $genus(E(K)) = g $ then $g = t(K) +1$ and assume also that  $cl(\partial_+ V_1 - A)$ has $k$ components 
$\Sigma_i$ each of which is a surface of genus $g_i$ with $c_i$ boundary components. Since we can obtain
$\Sigma_i$ from $\Sigma$ = $\partial_+ V_i$  by cutting $\Sigma$ along annuli we have the following formula:
\vskip10pt
\begin{itemize} 
\item[\rm (a)]$ -\chi(\Sigma) =  2g - 2 = 2t(K)$
 \item[\rm (b)]$2(2+ d) = \sum\nolimits_{i = 1}^k c_i$
\item[\rm (c)]$\sum\nolimits_{i = 1}^k c_i = 2t(K) + \sum\nolimits_{i = 1}^k 2 - 2g_i$ \hskip 8pt Thus:
\item[\rm (d)]$2(2+ d) = 2t(K) + \sum\nolimits_{i = 1}^k 2 - 2g_i$.
\end {itemize}
Given an $n_i$-float $F_i$ contained in $E(K_j)$ it determines $n_i - 1$ tunnels which intersect $A$
and an additional $g_i$ interior tunnels which do not meet $A$. Recall that since we have two components 
which intersect $A$ in two vertical annuli each the total number of floats is $k - 2$. Further note that each
disk in $A \cap V_1$ corresponds to tunnels on both sides of $A$. Putting all of the above together we see
that the total number of the tunnels in both induced Heegaard splittings on $E(K_1)$ and $E(K_2)$ is:

$$t^*(K_1) + t^*(K_2) = 2d -(k - 2) + \sum\nolimits_{i = 1}^k g_i $$

\noindent Where $t^*(K_i )$ is the tunnel number of the induced Heegaard splitting on $E(K_i)$.
Combined with the above we get:

$$ t^*(K_1) + t^*(K_2) =  2t(K)  + (k - 2) -  \sum\nolimits_{i = 1}^k g_i$$
\vskip 10pt 
In the next section  we show that by ruling out $2$-floats of genus zero we can give a global upper 
bound on the number $(k - 2) - \sum\nolimits_{i = 1}^k g_i  $, and hence a meaningful upper bound
on $t^*(K_1) + t^*(K_2)$ in terms of  $t(K)$ (see Theorem \ref{bound}).
\end{remark}
\vskip 10pt

\section{Ruling out genus zero  2-floats} 

\vskip 13pt

We have the following theorem:

\begin {theorem}

\label {UturnProp}

Let $K_i \subset S^3$ be prime knots and $K = K_1 \# K_2$ be 
their connected sum. Let $(V_1, V_2)$  be a Heegaard splitting of  $E(K)$ and $A$ a decomposing
annulus which intersects $V_1$ in a minimal number of disks. Then: $V_1 - N(A)$  has at most one 
component which is a $2$-float of genus zero. 

\end{theorem}

\noindent {\bf Remark:} If we generalize the definition of tunnels to a neighborhood of a complete system of 
meridian disks for $V_1$ we have: With the above assumptions at most one tunnel of $K$ has a  
U-turn. i.e., at most one tunnel pierces $A$ in one direction and then pierces it again in the
opposite direction without meeting $V_1$. Furthermore this phenomenon occurs at most once.

\vskip10pt

\begin {proof} Assume  that some 2-float is a $3$-ball meeting $A$ 
in exactly two disks. Denote these disks by $D_1$ and $D_2$ but note that these indices do not
necessarily agree with the natural order defined on the disks $D_i$ by  Theorem \ref {darcThm}.  We
first need the following lemma:

\vskip10pt

\begin {lemma}

\label{DiskLemma}
There is some $E_j \subset {\cal E}$ and a sub-disk $\Delta \subset E_j$ so that 
$\partial \Delta \subset 2-float  \cup A$. Hence $\partial \Delta = \cup (\alpha_r \cup \beta_s)$, where the
$\beta_s$ arcs are contained in the 2-float and the $\alpha_r $ are arcs on $A$. 
\end{lemma}

\begin {proof} (of Lemma) By Theorem \ref {darcThm} there is some arc $\alpha$ in some $E_j$ with 
end points on $D_1$.  Consider the disk $\Delta$ of $E_j - A$ adjacent to $\alpha$ on the same side of 
$A$ as the 2-float. We can assume that $\Delta$ is to the right of $\alpha$ on $E_j$ and it cannot be
outer-most  as the 2-float has no genus.  Hence there are more arcs of $E_j \cap A$ further to the right of
$\alpha$. If all  such arcs of $\Delta$ which are not on $\partial E_j$, are of type III the disk 
$\Delta$ satisfies the conclusion of the lemma and we are done: Since all arcs of type III have end
points on $D_1$ or $D_2$ the $\beta$ arcs must be contained in the 2-float.

So we assume that some arc to the right of $\alpha$ is of type I. If the lemma  fails there are at 
least two  arcs of type I with  one end point not on $D_1$ or $D_2$. Consider such an
arc $\rho$ of type I, it cannot be outer-most by Theorem \ref {darcThm}. Hence further to the right there is
some arc $\alpha'$ of type III with end  points on $D_1$ or $D_2$. One of the two disks adjacent to
$\alpha'$ is on the 2-float. Assume that it is $\Delta'$ and that it is on the right of $\alpha'$. We now start
the argument again with $\alpha'$. This procedure must end since the intersection is finite. 

Assume therefore that $\Delta'$ is to the left of $\alpha'$. If all arcs in $\partial \Delta' -\partial E_j$
have end points on $D_1$ or $D_2$ we are done as before. If there is an arc with no end points on 
$D_1$ or $D_2$ then there is an arc $\rho'$ of type I  with exactly one end point on $D_1$ or $D_2$.
It cannot be outermost as before so farther out there is an arc of type III with end points on
$D_1$ or $D_2$. So we can start the argument with  $\rho'$. However the procedure must
terminate as the intersection is finite. Hence at some stage we obtain a disk $\Delta$ with
$\partial \Delta - \partial E_j$ consisting of arcs of type I or type III all of which have 
end points on $D_1$ or $D_2$. Hence all the $\beta$ arcs are on the 2-float. (see Fig. 1)

\end {proof}

\begin {corollary} The disk $\Delta$ contains at most one arc $\alpha$ of type I on its boundary. 

\end {corollary}

\qed

\bigskip

Consider now an essential sub-annulus $A'$ of $A$ containing the disks $D_1$ and 
$D_2$. It is a meridional annulus in $(S^3, K)$ so we can cap off $A'$ by two meridian disks 
$D^*_1$ and $D^*_2$ in $(S^3, K)$ to obtain a 2-sphere intersecting $K$ in exactly two points in
$D^*_1$ and  $D^*_2$. If we attach the boundary of the 2-float to this 2-sphere along 
$D_1$ and $D_2$ we get a 2-torus $T$. By the above lemma $\partial \Delta$ is contained in $T$.
.
\bigskip

\noindent \underbar{Case 1:} Assume that $\partial \Delta$ is an inessential curve on $T$
and bounds a disk $\Delta'$ there which does not contain the disks $D^*_1$ and  $D^*_2$. 
The intersection of
$\partial \Delta$ with a core curve of the meridional annulus $A'$ is even. Similarly the intersection of $\partial
\Delta$  with the boundary of a  cocore disk of the 2-float is even. Hence the number of arcs of type I is
even and so is the number of $\beta$ arcs (these are the arcs which intersect the boundary of a  
cocore disk of the 2-float). Hence the number of arcs of type III (the $\alpha$ arcs) is also even.
As a consequence the disk   $\Delta'$ is a union of bands glued together to each other at their ends.
The bands correspond to the areas in $A$ between the arcs of type I and between the arcs of type III
and also on the 2-float between the $\beta$ arcs (see Fig. 4).

Since the bands are glued to each other along small arcs on both ends, the number of gluing arcs
is equal to the number of bands. So an Euler characteristic argument shows that

\bigskip
\centerline {$ \chi (\Delta') = \sum \chi(bands) -  \sum \chi (gluing \hskip3pt arcs) = 0$} 

\bigskip

\noindent  But this is obviously a contradiction and hence $\partial \Delta$ is essential in $T$  or bounds a 
disk $\Delta'$ containing one or both of $D^*_1$ and  $D^*_2$.

\begin{figure}[htm]
\hskip20pt{\epsfxsize=4truein\epsfbox{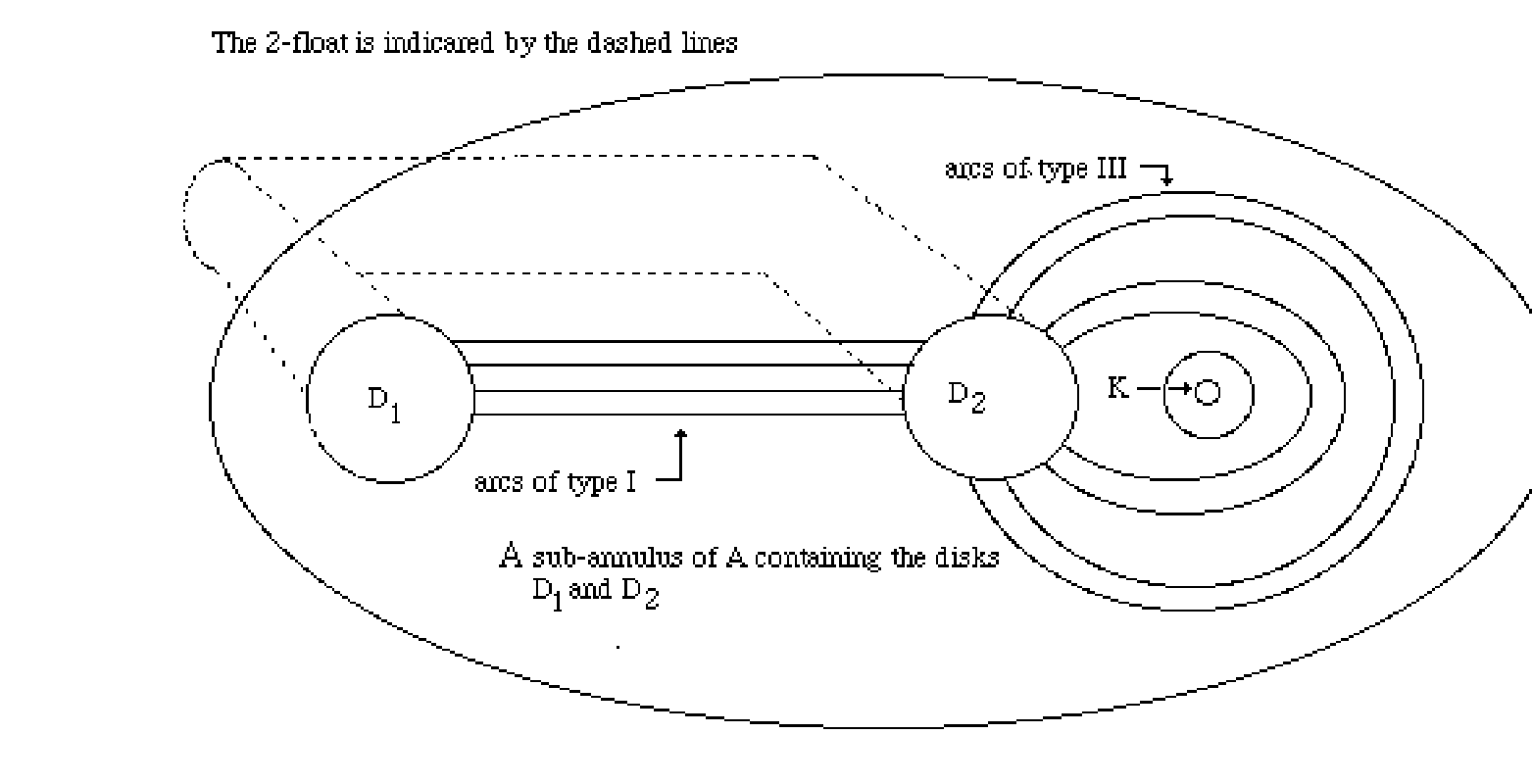}}

\centerline {Fig. 4}
\end{figure}
\vskip10pt 

\noindent \underbar{Case 2:}   If $\Delta'$ contains only one of $D^*_1$ or $D^*_2$ then the 
$2$-sphere $\Delta' \cup \Delta$ intersects $K$ in a single point in contradiction. Assume therefore,
that  $\Delta'$ contains both of $D^*_1$ and $D^*_2$. Consider the annulus 
$A'' = cl((\Delta' \cup \Delta) - (D^*_1 \cup D^*_2))$. We can push $A''$ by a small ambient 
isotopy relative to its boundary  $\partial D^*_1 \cup \partial D^*_2 $ in the direction away from 
the $2$-float so that after the isoptopy $A''$ does not intersect the disks $D_1$ and $D_2$. 
If  the sub-annuls $A - A'$ contains disks of intersection with $V_1$ then the annulus $A''$ is not parallel into
$A - A'$ and the annulus $(A -A') \cup A''$  is a new  decomposing annulus intersecting $V_1$ in fewer 
disks than $A$, in contradiction to the  choice of $A$.
\vskip5pt

\noindent \underbar{Case 3:}
If $\partial \Delta$ is an essential curve on $T$ we do 2-surgery on $\partial \Delta$ by removing 
an annulus neighborhood of $\partial \Delta$  and gluing two copies of $ \Delta$. By an Euler 
characteristic argument we obtain a 2-sphere intersecting $K$ in two points on  $D^*_1$ and 
$D^*_2$. If we remove  $D^*_1$ and $D^*_2$ we obtain an annulus $A''$. As above can now 
replace the annulus $A'$ by the annulus $A''$ and get a new  annulus $(A - A') \cup A''$ which 
does not intersect the disks $D_1$ and $D_2$. The annulus $A''$ cannot be parallel into $A'$
as this would imply that $(A - A') \cup A''$ is parallel into $A$ in contradiction to the choice of $A$
as minimizing the intersection with $V_1$. If  the  sub-annuls $A - A'$ contains disks of intersection 
with $V_1$ then the annulus $A''$ is not  parallel into $A - A'$ and the annulus $(A -A') \cup A''$  is a new 
decomposing annulus  intersecting $V_1$ in fewer  disks than $A$, in contradiction to the  choice of $A$. 

Hence we conclude that  the disks $D_1$ and $D_2$ must be the first and last disks in $A \cap V_1$,
i.e., the two disks closest to $A_1^*$ and $A_2^*$. This implies that each knot complement $E(K_1)$
and $E(K_2)$ contains at most one component of $V_1 - A$ which is a $2$-float of genus zero. Since
$V_1$ is connected we cannot have two $2$-floats of genus zero one on each side both intersecting
$A$ in the disks $D_1$ and $D_2$.

\end {proof}

As an immediate application of Theorem \ref{UturnProp} we have the following bound on the tunnel
number of $K = K_1 \# K_2$ in terms of the tunnel numbers of $K_1$ and $K_2$.  Note that
in [SS] Schultens and Scharlemann have obtained a better bound. However we bring it here as it  is an
immediate corollary of our main theorem:

By the above theorem at most one  component of $V_1 - A$ is a $2$-float. Hence for all other components
$C_i$ we have $\partial C_i - A = \Sigma_i$ is a surface so that $-\chi(\Sigma_i) \geq 1$, as $\Sigma_i$ is
either a punctured $2$-sphere with at least three punctures or a surface of positive genus with at least one
puncture. Since we also have, by Proposition \ref{innertunnelProp}, at least one component with positive
genus we have as a worst case the following situation:

One component $C_1$ of $V_1 - A$ which intersects the vertical annuli and $\partial C_1 - \partial
E(K)$ is an  annulus. Since $V_1$ is connected the other component  $C_1$ of $V_1 - A$ which intersects the
vertical annuli must have at least one other disk of intersection with $A$. Furthermore we have one component
which is a $2$-float, one component of genus one and all other components are $3$-floats of genus zero. For a
closed surface of genus $h$ the latter type of components can number at most $2h - 2$ so $k = 2(g - 2)  - 2 + 4 =
2g - 2$ and
 $\sum\nolimits_{i = 1}^k g_i  -(k - 2)  \geq 1- (2g - 2) +2 = 3 - 2t(K)$. Hence by Remark \ref{benefit}:

\begin{theorem}\label{bound} 

For all prime knots  $K_1$ and $K_2$ in $S^3$ and $K = K_1 \# K_2$ we have: 

$$t^*(K_1) + t^*(K_2)\leq 4 t(K)  - 3 $$  \hskip80pt \qed
\end{theorem}

The above theorem suggest that by ruling out $n$-floats of genus zero with $n = 3,4. \dots$ one can get better
bounds  which would converge to approximately $t^*(K_1) + t^*(K_2)\leq 2 t(K) $.

\eject

\section {References}
\bigskip

\noindent [Ja] \hskip 28pt W. Jaco; {\it Lectures on three manifold topology} CBMS Regional 
\vskip8pt

\noindent \hskip 49pt Conference  Series   in Mathematics 43 1977.
\vskip8pt

\noindent [Ko] \hskip24pt T. Kobayashi; {\it A construction of arbitrarily high degeneration 

\vskip8pt

\hskip30pt of tunnel number of knots under connected sum}, J. Knot Theory 

\vskip8pt

\hskip30pt and its Ramifications 3 (1994), 179 - 186 .

\vskip8pt

\noindent [Kw] \hskip21pt H.-Z. Kwong; {\it Straightening tori in Heegaard spkittings}, Ph.D 

\vskip 8pt

\hskip 30pt Thesis U.C. Santa Barbara 1994 . 

\vskip8pt

\noindent [MR]\hskip 19pt  Y. Moriah, H. Rubinstein; {\it Heegaard structures of negatively curved

\vskip8pt

\noindent \hskip 46pt  3-manifolds} Comm. in Anal. and Geom. 5 (1997), 375 - 412 .

\vskip8pt

\noindent [Mo1]\hskip 20pt K. Morimoto; {\it  On the additivity of tunnel number of knots}, Topol- 

\vskip8pt

\noindent \hskip 50pt ogy and its Applications 53 (1993)  37 - 66.

\vskip8pt

\noindent [Mo2]\hskip 24pt K. Morimoto; {\it  There are knots whose tunnel numbers go down un- 

\vskip8pt

\noindent \hskip 48pt der connected sum}, Proc. Amer. Math. Soc. 123 (1995) 3527 - 

\vskip8pt

\noindent \hskip 50pt 3532.

\vskip8pt

\noindent [Mo3]\hskip 24pt K. Morimoto; {\it  On the super additivity of tunnel number of knots},

\vskip8pt

\noindent \hskip 52pt preprint

\vskip8pt

\noindent [MSY] \hskip 10pt K. Morimoto, M. Sakuma, Y. Yokota; {\it  Examples of tunnel number 

\vskip8pt

\hskip 28pt  one knots which have the property that \lq\lq 1 + 1 = 3"},  Math. Proc. 

\vskip8pt

\noindent \hskip 45pt   Camb. Phil. Soc., 119 (1996) 113 - 118 .

\vskip8pt

\noindent [Oc]\hskip 28pt  M. Ochiai; {\it On Haken's theorem and its extension}, Osaka J. of

\vskip8pt

\noindent \hskip 48pt Math. 20 (1983) 461 - 468.

\vskip8pt

\noindent [SS] \hskip 24pt M. Scharlemann, J. Schultens; {\it Annuli in generalized Heegaard 

\vskip8pt

\hskip30pt splitting and degeneration of tunnel number}, Preprint.

\vskip35pt

\obeyspaces   Yoav Moriah

\obeyspaces   Department of Mathematics

 \obeyspaces  Technion, Haifa  32000,

\obeyspaces  Israel

\vskip10pt

 ymoriah@tx.technion.ac.il

\vskip20pt

\end{document}